\documentclass[conference]{IEEEtran}
\IEEEoverridecommandlockouts
\usepackage{cite}
\usepackage{amsmath,amssymb,amsfonts}
\usepackage{algorithm2e}
\usepackage{hyperref}
\usepackage{graphicx}
\usepackage{subcaption}
\usepackage{textcomp}
\usepackage{xcolor}
\usepackage{bm}
\def\BibTeX{{\rm B\kern-.05em{\sc i\kern-.025em b}\kern-.08em
    T\kern-.1667em\lower.7ex\hbox{E}\kern-.125emX}}

\newcommand{\bA}{\mathbf{A}}
\newcommand{\bD}{\mathbf{D}}
\newcommand{\bb}{\mathbf{b}}
\newcommand{\bc}{\mathbf{c}}
\newcommand{\bd}{\mathbf{d}}
\newcommand{\br}{\mathbf{r}}
\newcommand{\bw}{\mathbf{w}}
\newcommand{\bx}{\mathbf{x}}

\begin{document}

\title{Stable Iterative Solvers \\for Ill-conditioned Linear Systems
}

\author{\IEEEauthorblockN{Vasileios Kalantzis, Mark S.\ Squillante, Chai Wah Wu}
\IEEEauthorblockA{\textit{Mathematics of Computation} \\
\textit{IBM Research} \\
\textit{Thomas J.\ Watson Research Center, 
Yorktown Heights, NY, USA} \\
vkal@ibm.com, mss@us.ibm.com, cwwu@us.ibm.com}
}

\maketitle

\begin{abstract}
Iterative solvers for large-scale linear systems such as Krylov subspace methods can  diverge when the linear system is ill-conditioned, thus significantly reducing the applicability of these iterative methods in practice for high-performance computing solutions of such large-scale linear systems. To address this fundamental problem, we propose general algorithmic frameworks to modify Krylov subspace iterative solution methods which ensure that the algorithms are stable and do not diverge. 
We then apply our general frameworks to current implementations of the corresponding iterative methods in \verb+SciPy+ and demonstrate the efficacy of our stable iterative approach with respect to numerical experiments across a wide range of synthetic and real-world ill-conditioned linear systems. 
\end{abstract}


\section{Introduction}
Consider the computationally intensive problem of solving large systems of linear equations, namely given a matrix $\bA$ and a vector $\bb$, find the vector $\bx$ 
such that $\bA\bx=\bb$ or such that $\|\bA\bx-\bb\|$ is minimized.
The solution of this general class of problems is an operational workhorse in the 
high-performance computations of modern science and engineering. 
Iterative solvers such as Krylov subspace methods have emerged as an important 
general class of algorithms for solving these large-scale linear system problems, especially in the case of sparse linear systems~\cite{Saad2003,liesen2013krylov}. 
However, when the linear system is ill-conditioned, i.e., the condition number of the matrix $\bA$ is large, such Krylov subspace iterative methods can often diverge which leads to numerical overflow or inaccurate results. This can happen even when preconditioners are used, if the preconditioner and $\bA^{-1}$ are mismatched.
Examples of ill-conditioned linear systems found in practice include the numerical solution of the high-frequency Helmholtz equation in 3D~\cite{xi2017rational} and the shifted matrices in rational filtering preconditioners~\cite{austin2024rational}.
In addition, various examples of ill-conditioned matrices are also prevalent in a standard library of matrices for real-world linear systems~\cite{sparse-ufl2011}.

Since ill-conditioned linear systems arise naturally in practice and their occurrences are not uncommon, the issues of inaccurate results, or even divergence, that can be exhibited by classical Krylov subspace iterative solvers when applied to such large-scale ill-conditioned systems represents a fundamental problem in high-performance computing.
This in turn can significantly reduce the applicability of classical Krylov subspace iterative methods in practice for high-performance solutions of these large-scale ill-conditioned systems. 
Hence, our goal is to address the fundamental problems associated with solving large-scale ill-conditioned linear systems using Krylov subspace iterative methods.
Examples of relevant Krylov subspace methods of interest include various forms of 
the
conjugate gradient (CG)~\cite{cg1950,cg1952}
algorithm
and the
generalized minimal residual (GMRES)~\cite{gmres1986,Saad2003} algorithm.

In this paper, we aim to mitigate the above fundamental problems by introducing general algorithmic frameworks to augment the broad class of Krylov subspace iterative methods in a manner that solves large-scale ill-conditioned linear systems while ensuring stability and guaranteeing a lack of divergence.
Beyond our algorithmic contributions, we conduct extensive numerical experiments that demonstrate the significant benefits of our general frameworks with respect to various forms of CG and GMRES and related methods.
In particular, we first show how classical Krylov subspace and related iterative methods diverge when the linear system is sufficiently ill-conditioned, further quantifying these effects.
We then demonstrate that the corresponding versions of these iterative methods augmented with our general algorithmic frameworks are stable, do not diverge, and provide significantly more accurate results.
Moreover, for less ill-conditioned linear system such that the classical Krylov iterative method does not diverge, our numerical experiments show that the Krylov subspace iterative methods augmented with our general algorithmic frameworks often provide significantly better accuracy than the original Krylov iterative method and never perform worse than the original method.

We note that Krylov subspace iterative methods have also found applications in high-performance computing environments that are based on inaccurate computations, such as those employing analog (see, e.g.,~\cite{HPEC2023}) and mixed-precision (see, e.g.,~\cite{Lindquist:2021}) technologies.
Although our general algorithmic frameworks introduced herein for Krylov subspace iterative methods can also be exploited to mitigate the problems of ill-conditioned linear systems in such inaccurate computing environments, we do not directly consider these computing environments in this paper.

The remainder of this paper is organized as follows.
We first present in Section~\ref{sec:frameworks} our general algorithmic frameworks.
We then provide in Section~\ref{sec:experiments} representative samples from extensive numerical experiments across a wide range of synthetic and real-world ill-conditioned linear systems based on various classical Krylov subspace
iterative linear system solvers, followed by concluding remarks in Section~\ref{sec:conclusions}.

\section{Our General Algorithmic Frameworks}
\label{sec:frameworks}
Before presenting our general stable iterative algorithmic frameworks for large-scale ill-conditioned linear systems, it is important to begin by considering classical Krylov subspace iterative solvers that, after starting with an initial guess $\bx_0$, operate by determining a vector $\bd_i$ at the $i$-th iteration to update the current solution $\bx_{i+1} = \bx_i + \bd_i$.
The residual of the solution at the $i$-th iteration is then defined as $\br_i = \bb-\bA\bx_i$. 
This classical iterative linear system solver framework is generically shown in Algorithm~\ref{alg:one},
the main steps of which can lead this classical iterative method to diverge when the matrix $\bA$ is ill-conditioned.

\RestyleAlgo{ruled}
\begin{algorithm}
\caption{Classical iterative linear system solver}
\label{alg:one}
\KwData{$\bA$, $\bb$, $\bx_0$}
\KwResult{$\bx$ such that $\bA\bx=\bb$}
\While{Stopping criteria is not satisfied}{
determine vector $\bd_i$\;
update $\bx_{i+1} \gets \bx_i + \bd_i$\;
compute residual $\br_{i+1} \gets \bb-\bA\bx_{i+1}$\;
}
\end{algorithm}

To address the divergence issues and related numerical inaccuracies of the classical Krylov iterative methods generically depicted in Algorithm~\ref{alg:one} when applied to ill-conditioned linear systems, our first general approach consists of augmenting this broad collection of Krylov subspace iterative linear system solvers with an algorithmic framework based on the inclusion of a line search along the direction $\bd_i$ in a manner that ensures the residual is reduced, i.e., ensuring $\|\bA\bx_{i+1}-\bb\|_2 \leq \|\bA\bx_i-\bb\|_2$, and that guarantees divergence does not occur under ill-conditioned linear systems.
More precisely, we modify the update of the solution $\bx_{i+1}$ at the $i$-th iteration to be $\bx_{i+1} = \bx_i + \alpha_i \bd_i$ which includes the additional scalar computed as $\alpha_i = \br_i^\top \bw_i/\|\bw_i\|_2^2$, where $\bw_i = \bA\bd_i$.
It is readily apparent that the desired inequality 
\begin{equation}\|\bA\bx_{i+1}-\bb\|_2 \leq \|\bA\bx_i-\bb\|_2
\label{eqn:residual}
\end{equation}
is
satisfied.
In addition, we further modify the update of the residual $\br_{i+1}$ at the $i$-th iteration to either be $\br_{i+1} = \br_i - \alpha_i \bw_i$, given that $\bw_i$ is already computed, or the residual can be occasionally recomputed via $\br_{i+1}=\bb-\bA\bx_{i+1}$, depending on the iterative method of interest. 
Our first general stable iterative framework is shown in Algorithm~\ref{alg:two}, where the modified main steps of our first algorithmic framework ensure a stable iterative method that does not diverge when the matrix $\bA$ is ill-conditioned.

\RestyleAlgo{ruled}
\begin{algorithm}
\caption{First stable iterative linear system framework}
\label{alg:two}
\KwData{$\bA$, $\bb$, $\bx_0$}
\KwResult{$\bx$ such that $\bA\bx=\bb$}
\While{Stopping criteria is not satisfied}{
determine vector $\bd_i$\;
compute $\bw_i = \bA\bd_i$ and $\alpha_i = \br_i^\top \bw_i/\|\bw_i\|_2^2$\;
update $\bx_{i+1} \gets \bx_i + \alpha\bd_i$\;
update residual $\br_{i+1} \gets \br_i - \alpha\bw_i$ or compute residual $\br_{i+1} \gets \bb-\bA\bx_{i+1}$\;
}
\end{algorithm}

We next extend our first general stable iterative algorithmic framework to consider augmenting the broad class of Krylov subspace iterative methods with the inclusion of a line search along the directions of both $\bx_i$ and $\bd_i$ in a manner that minimizes the residual-norm $\|\bA\bx_{i+1}-\bb\|_2$ while continuing to guarantee that Eq.~(\ref{eqn:residual}) is satisfied and divergence does not occur under ill-conditioned linear systems.
More precisely, we first construct the $n \times 2$ matrix $\bD_i = [\bx_i;\bd_i]$.
Then, since $\bD_i^\top \bA^\top \bA \bD_i$ can be singular or nearly singular, we obtain the least square solution $\bc_i$ of the $2\times 2$ system
\begin{equation}\label{eq:extended-framework}
    \bD_i^\top \bA^\top \bA \bD_i\bc_i = \bD_i^\top\bA^\top\bb,
\end{equation}
rather than solving it exactly, and update accordingly the solution $\bx_{i+1}$ and the residual $\br_{i+1}$ at the $i$-th iteration.
Our second extended general stable iterative framework is shown in Algorithm~\ref{alg:three}, where the modified main steps of our second algorithmic framework further ensure a stable iterative method that does not diverge when the matrix $\bA$ is ill-conditioned.

\RestyleAlgo{ruled}
\begin{algorithm}
\caption{Second stable iterative linear system framework}
\label{alg:three}
\KwData{$\bA$, $\bb$, $\bx_0$}
\KwResult{$\bx$ such that $\bA\bx=\bb$}
\While{Stopping criteria is not satisfied}{
determine vector $\bd_i$\;
construct $n \times 2$ matrix $\bD_i = [\bx_i;\bd_i]$\;
obtain least square solution $\bc_i$ of $\bD_i^\top \bA^\top \bA \bD_i\bc_i = \bD_i^\top\bA^\top\bb$\; 
update $\bx_{i+1} \gets \bD_i \bc_i$\;
update $\br_{i+1} \gets \bb-\bA\bx_{i+1}$\;
}
\end{algorithm}

\section{Numerical Experiments}
\label{sec:experiments}
We now present an extensive collection of numerical experiments that apply our general stable algorithmic frameworks within the context of several of the classical Krylov subspace iterative linear system solvers found in \verb+SciPy+ \cite{2020SciPy-NMeth},
which
is the facto standard module for scientific computing in \verb+Python+ and contains various forms of Krylov subspace iterative methods such as 
CG~\cite{cg1950,cg1952}, 
biconjugate gradient (BICG)~\cite{Saad2003,liesen2013krylov}, 
biconjugate gradient stabilized (BICGSTAB)~\cite{vorst:bicgstab:1992,Saad2003,liesen2013krylov}, 
conjugate gradient squared (CGS)~\cite{Saad2003,liesen2013krylov}, 
GMRES~\cite{gmres1986,Saad2003}, 
loose GMRES (LGMRES)~\cite{ding2007loose,liesen2013krylov},
transpose-free quasi minimal residual (TFQMR)~\cite{freund:tfqmr:1993,liesen2013krylov},
and so on.
We therefore modify these iterative sparse linear system solvers in \verb+SciPy+ version 1.15.1 based on our general algorithmic frameworks to ensure that these iterative methods do not diverge.
Specifically, we implemented according to Algorithm~\ref{alg:two} and Algorithm~\ref{alg:three} stable versions of such iterative solvers in \verb+SciPy+ version 1.15.1, namely the algorithms
BICG (\verb+bicg+),
BICGSTAB (\verb+bicgstab+), 
CG (\verb+cg+),
CGS (\verb+cgs+),
GMRES (\verb+gmres+),
TFQMR (\verb+tfqmr+), and
LGMRES (\verb+lgmres+).

To demonstrate the efficacy of our general stable iterative frameworks, we conduct a large collection of numerical experiments that compare the solution-norm $\|\bx\|_2$ and the residual-norm $\|\bA\bx - \bb\|_2$ from the foregoing versions of the SciPy algorithms modified within the context of our general frameworks against those from the original algorithms in \verb+SciPy+.
This extensive collection of numerical experiments are based on three different sets of ill-conditioned matrices, that is instances of the class of Hilbert matrices, various ill-conditioned random matrices, and various ill-conditioned matrices taken from a standard library of matrices for real-world  systems~\cite{sparse-ufl2011}.
A representative sample of these numerical experiments are presented across the three sets of ill-conditioned matrices.

\subsection{Hilbert matrices}\label{sec:Hilbert}
\begin{figure}[htbp]
\centering
{\includegraphics[width=0.5\textwidth]{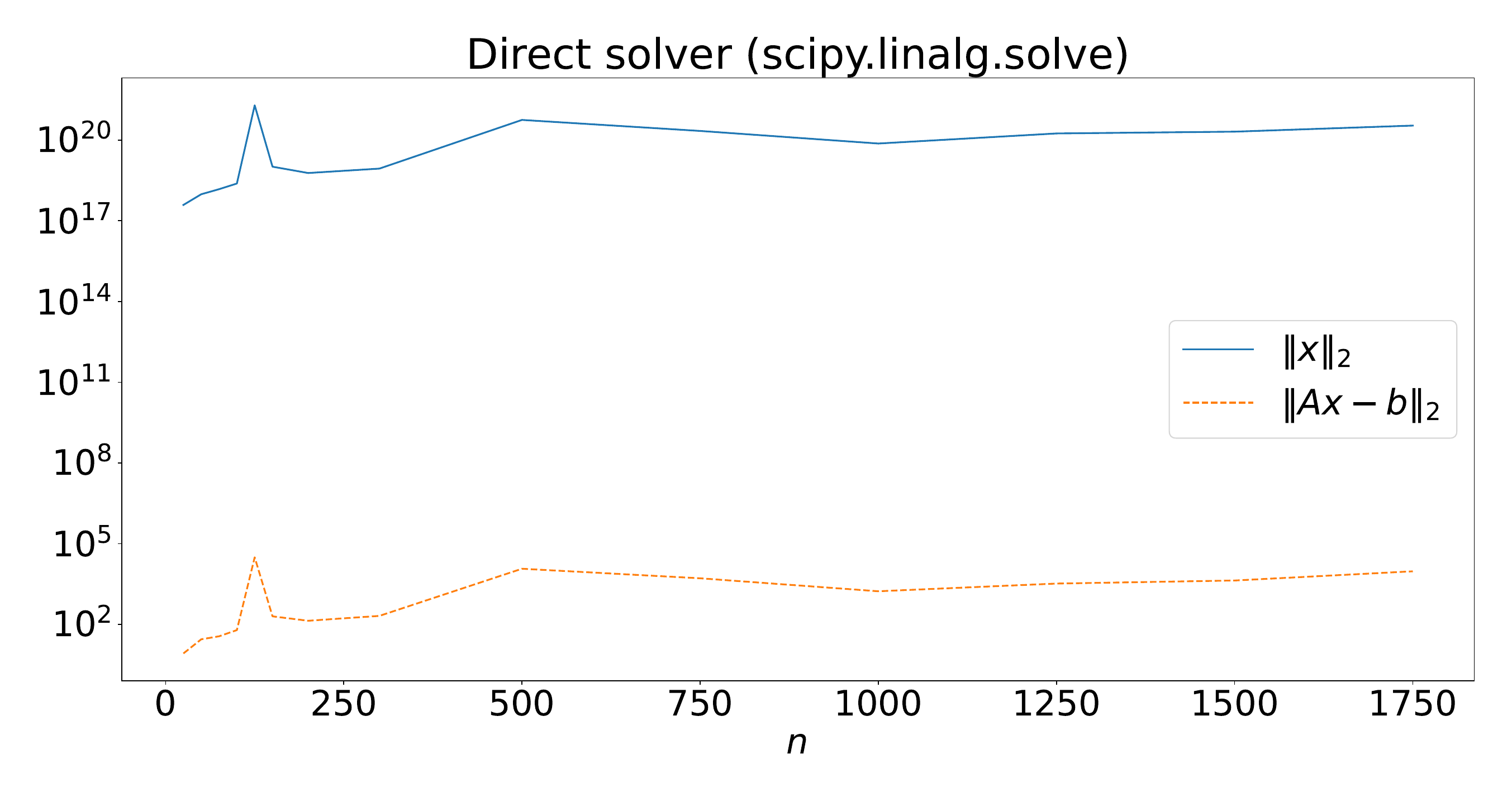}}
\caption{Solving linear systems with a Hilbert matrix using direct solver {\tt solve}.}
\label{fig:solve}
\end{figure}

Hilbert matrices form a well-known class of ill-conditioned matrices~\cite{Todd1954,Todd1961} that comprises $n\times n$ symmetric matrices defined as 
$$\bA_n(j,k) = \frac{1}{j+k-1}$$
for
$j,k\in [n] := \{1,\ldots,n\}$.
The inverses of the matrices $\bA_n$ are integer-valued and the condition numbers of $\bA_n$ grow as $O(\mu^n)$ for some constant $\mu \approx 33.97$, thus increasing very quickly with the dimension $n$ of the matrices $\bA_n$.
We randomly choose the vectors $\bb$ over $10$ trials (unless otherwise noted) and take the average of the Euclidean norm of the resulting solutions and residuals.
In Fig.~\ref{fig:solve}, we present numerical experiments which demonstrate that the direct solver SOLVE in \verb+SciPy+ results in relatively large values for both the solution-norm $\|\bx\|_2$ and the residual-norm $\|\bA\bx - \bb\|_2$ when $n$ increases
beyond relatively small values.

\begin{figure}[htbp]
\centering
{\includegraphics[width=0.5\textwidth]{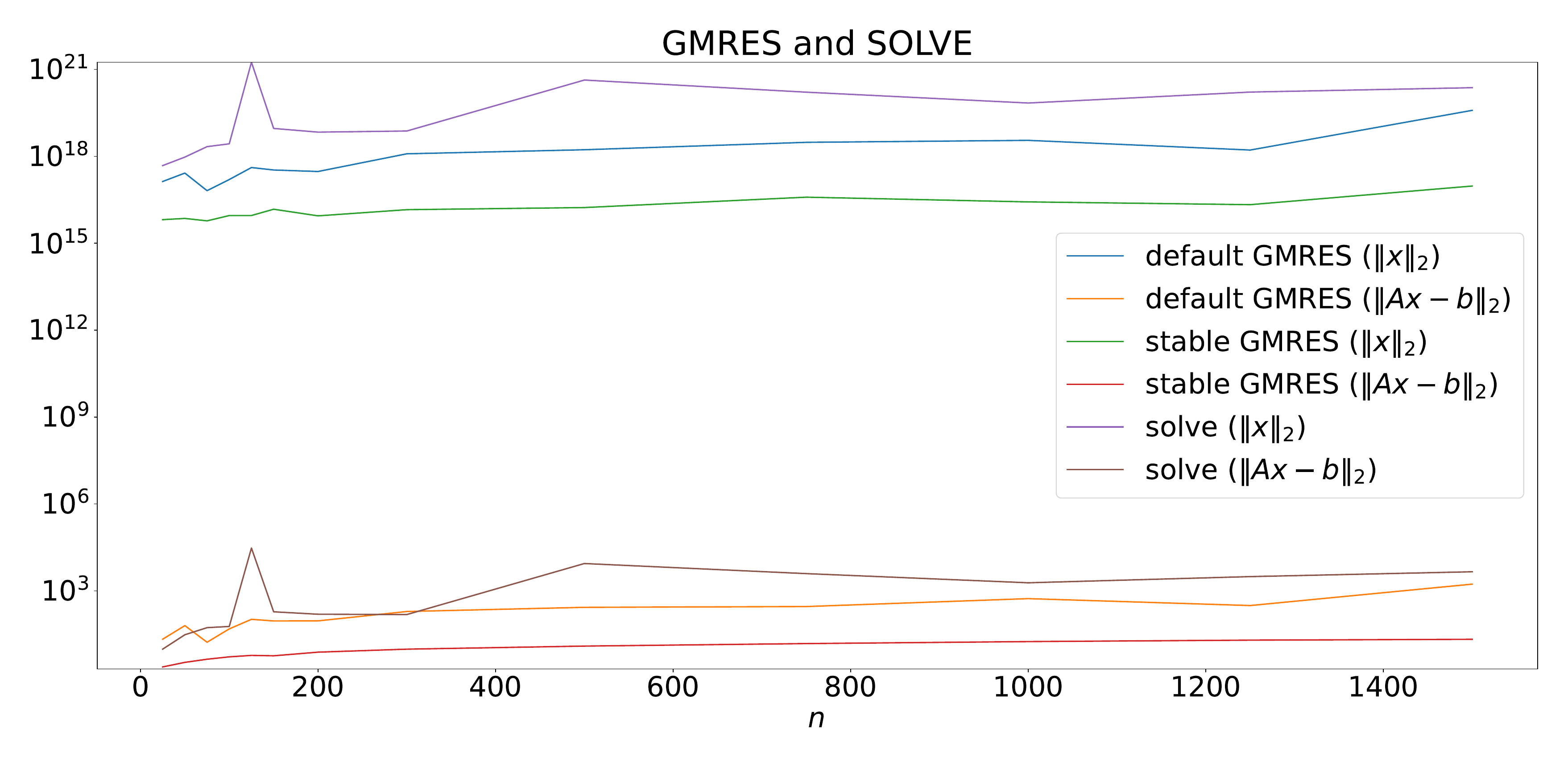}}
\caption{Solving linear systems with a Hilbert matrix using GMRES.}
\label{fig:gmres}
\end{figure}

We next show in Fig.~\ref{fig:gmres} the corresponding numerical results for the default GMRES algorithm and the stable GMRES algorithm modified according to our first framework, which illustrate that our stable algorithm ensures the solution norm and residual norm do not diverge for large $n$. 
For a firsthand comparison, the results for the direct solver SOLVE from Fig.~\ref{fig:solve} are also included in Fig.~\ref{fig:gmres}.
We observe that the numerical results in Fig.~\ref{fig:gmres} further show that our stable GMRES algorithm yields significantly smaller solution and residual norms compared with both the SOLVE and default GMRES algorithms.

\begin{figure}[htbp]
\centering
{\includegraphics[width=0.5\textwidth]{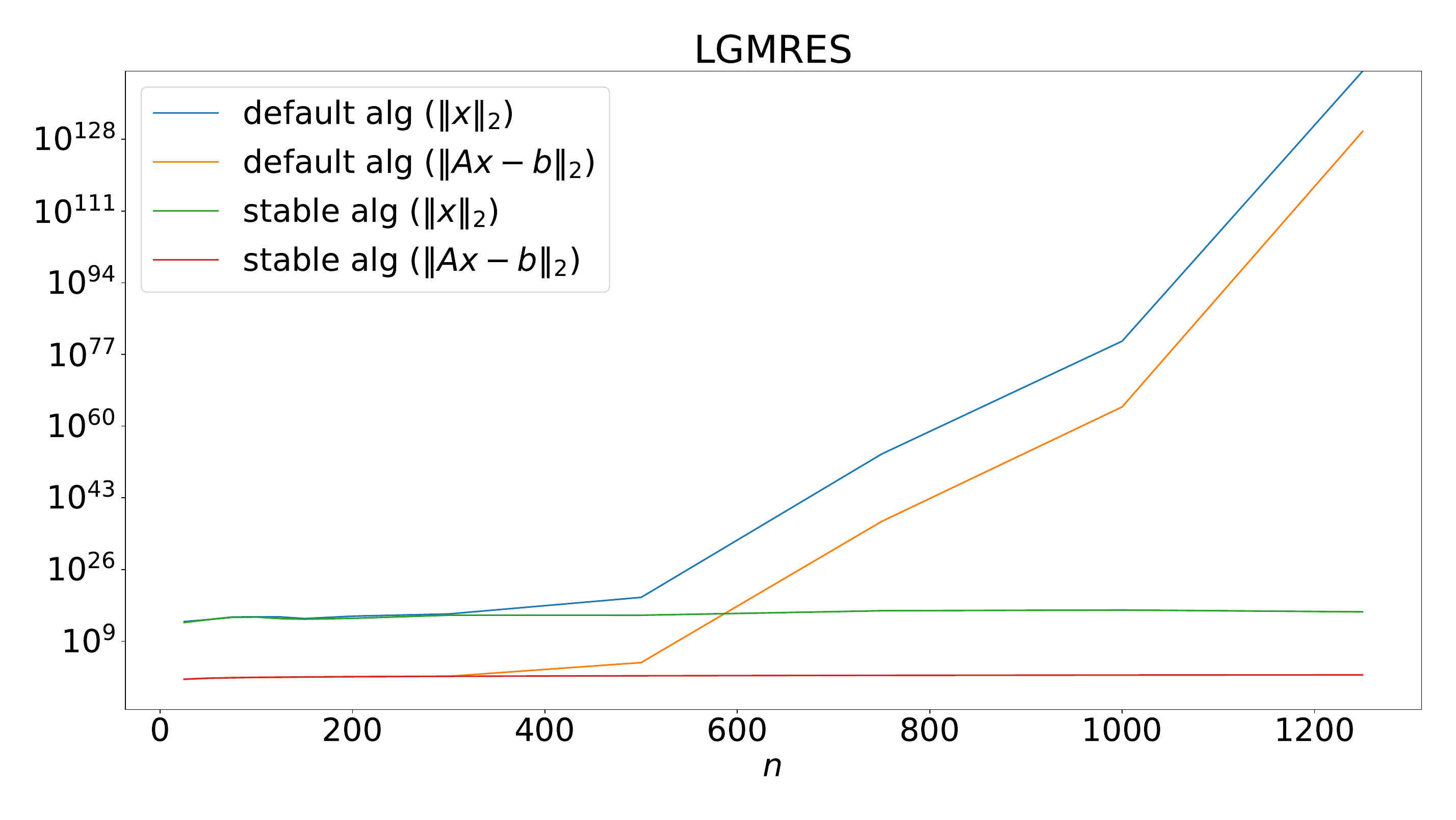}}
\caption{Solving linear systems with a Hilbert matrix using LGMRES.}
\label{fig:lgmres}
\end{figure}

\begin{figure}[htbp]
\centering
{\includegraphics[width=0.5\textwidth]{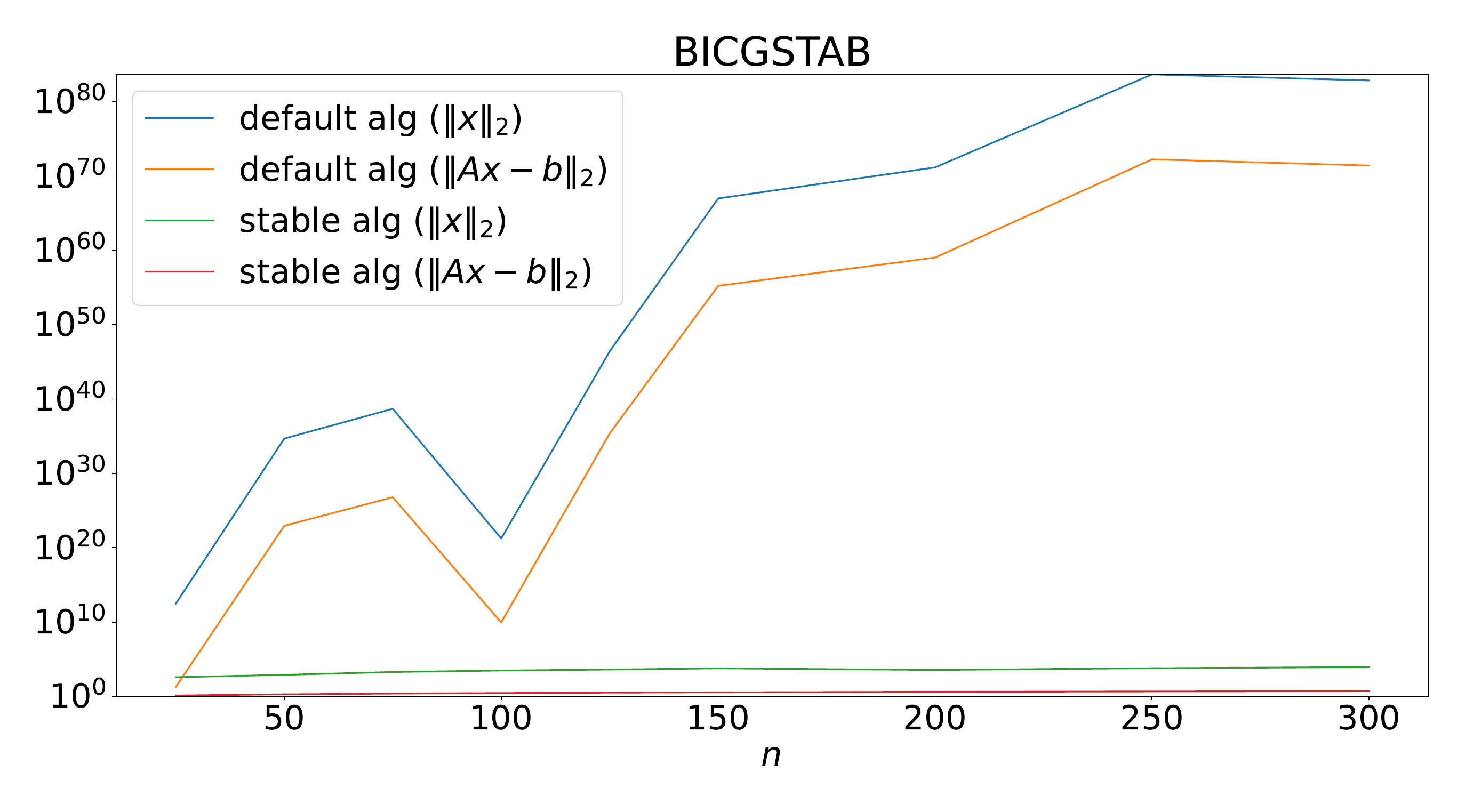}}
\caption{Solving linear systems with a Hilbert matrix using BICGSTAB.}
\label{fig:bicgstab}
\end{figure}

\begin{figure}[htbp]
\centering
{\includegraphics[width=0.5\textwidth]{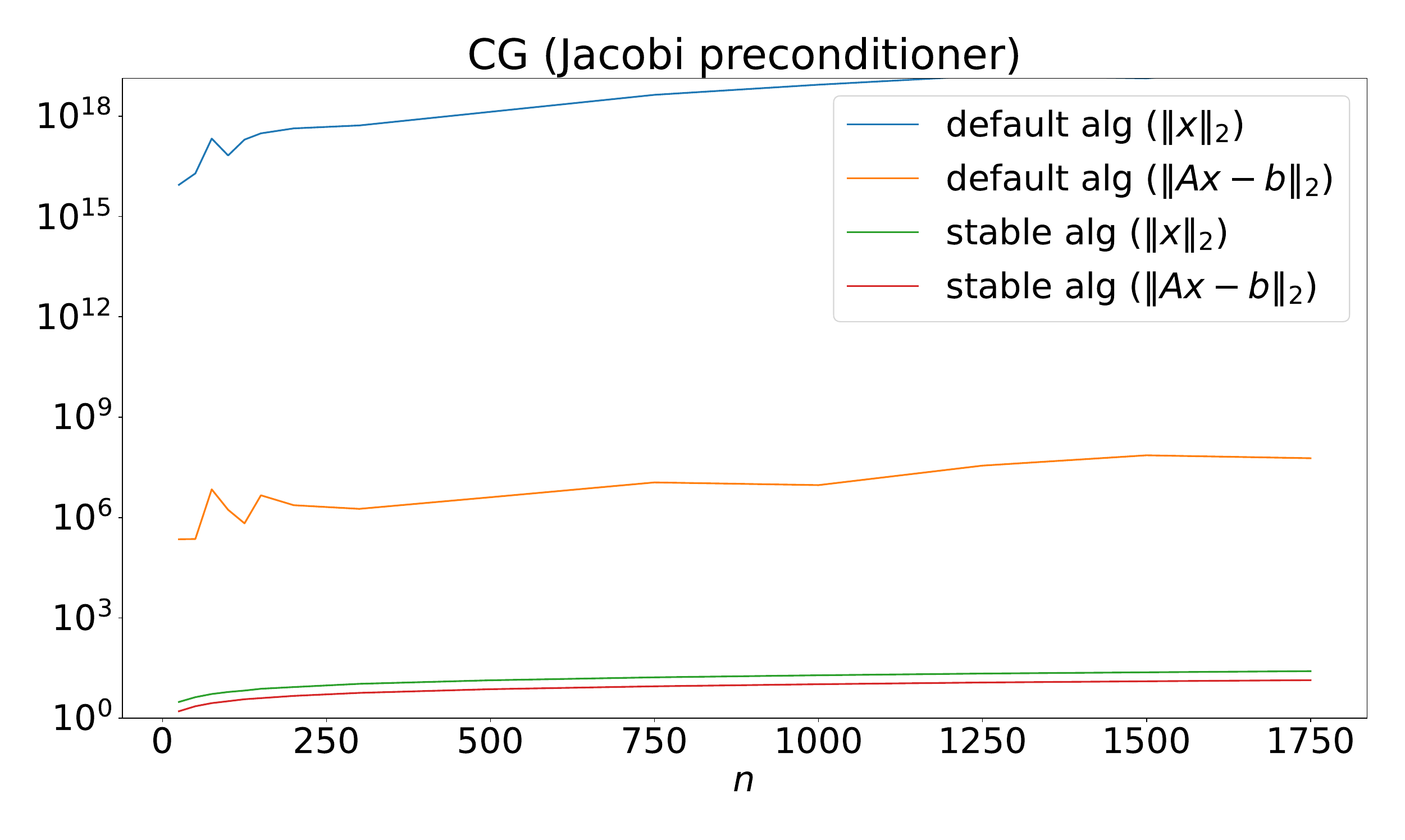}}
\caption{Solving linear systems with a Hilbert matrix using CG and the Jacobi preconditioner over $20$ trials.}
\label{fig:cg}
\end{figure}

\begin{figure}[htbp]
\centering
{\includegraphics[width=0.5\textwidth]{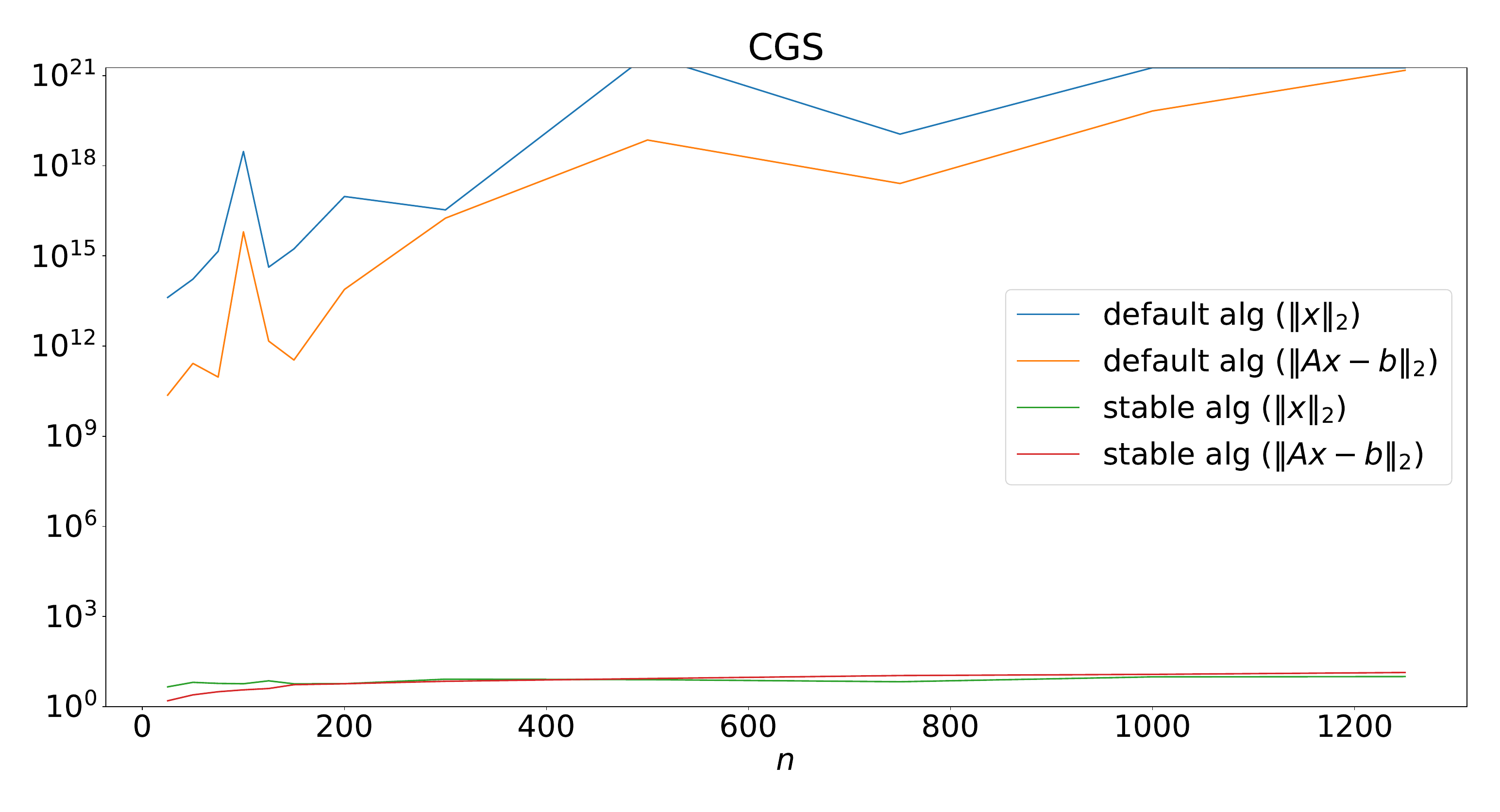}}
\caption{Solving linear systems with a Hilbert matrix using CGS.}
\label{fig:cgs}
\end{figure}

In a similar manner, we respectively show in Figs.~\ref{fig:lgmres}~--~\ref{fig:cgs} the results from the corresponding numerical experiments for the stable LGMRES, BICGSTAB, CG and CGS algorithms modified according to our first framework together with their default versions in \verb+SciPy+.
Furthermore, to show that the proposed frameworks are useful in preconditioned iterative solvers, the results for CG in Fig.~\ref{fig:cg} are obtained using the Jacobi preconditioner.
We observe from these numerical results that, in every case, the default algorithms diverge for large $n$, whereas our stable versions of the corresponding algorithms ensure the solution norm and residual norm do not diverge, remaining bounded and consistent across all values of $n$.
Moreover, we observe that our stable versions of the algorithms persistently provide significantly smaller solution and residual norms compared with the corresponding default algorithms.

\begin{figure}[htbp]
\centering
{\includegraphics[width=0.5\textwidth]{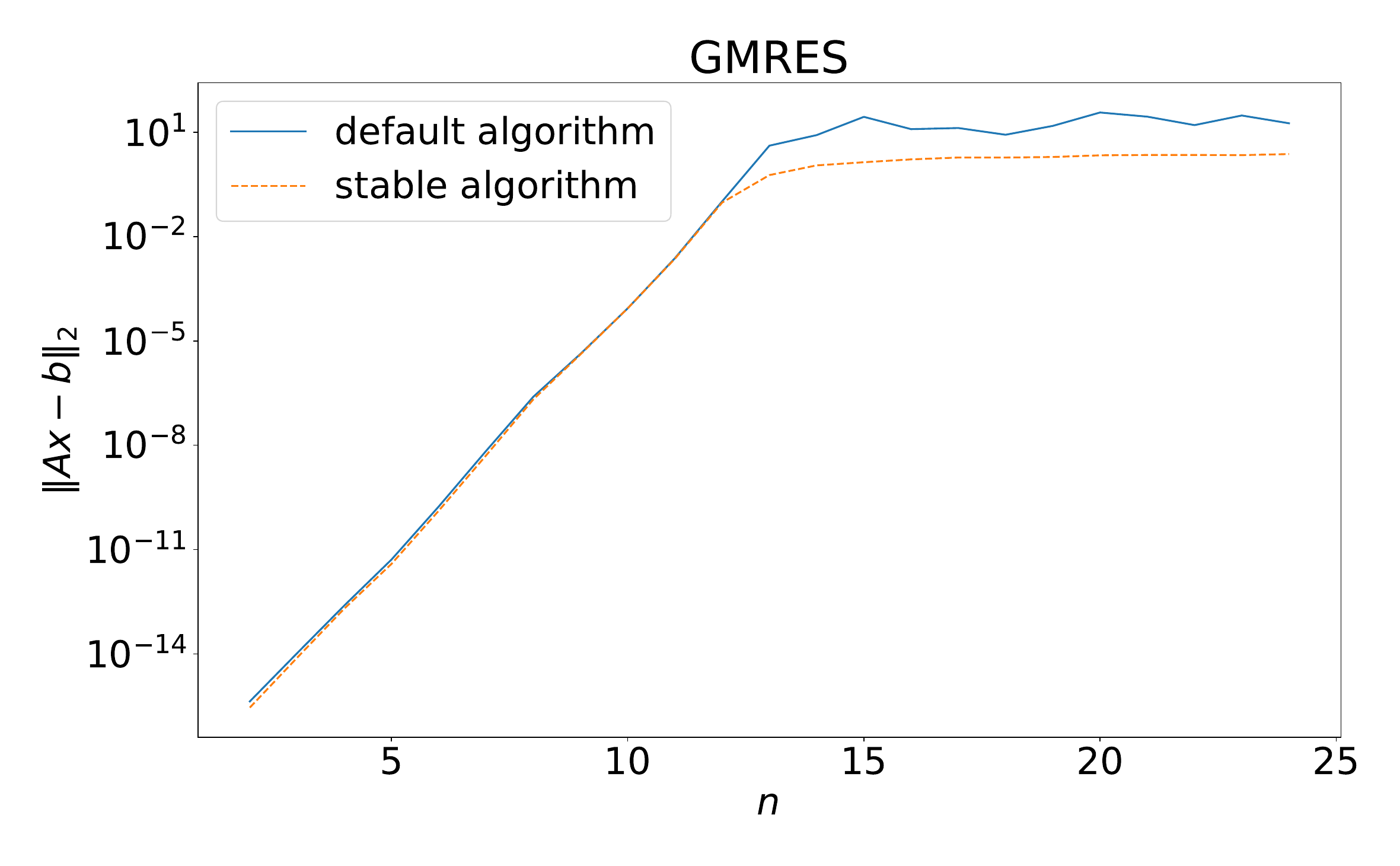}}
\caption{Solving linear systems with a Hilbert matrix using GMRES over 50 trials.}
\label{fig:gmressmalln}
\end{figure}

Zooming in on the numerical results from Fig.~\ref{fig:gmres} for small $n$, as a representative example of the illustrative behaviors from among the numerical experiments in Figs.~\ref{fig:gmres}~--~\ref{fig:cgs}, we demonstrate in Fig.~\ref{fig:gmressmalln} the stable behavior of the GMRES algorithm modified according to our first framework in comparison with its default counterpart in \verb+SciPy+.
We observe that, when $n$ is small, both the default GMRES algorithm and our stable GMRES algorithm perform equally well.
However, as $n$ increases, the default GMRES algorithm starts to perform considerably worse with larger and fluctuating residual norms compared to our modified GMRES algorithm which performs considerably better with smaller and stable residual norms.
As $n$ continues to increase, Fig.~\ref{fig:gmres} further shows the performance benefits of our stable GMRES algorithm, which consistently performs significantly better and in a more stable manner than the default GMRES algorithm.

\begin{figure}[htbp]
\centering
{\includegraphics[width=0.5\textwidth]{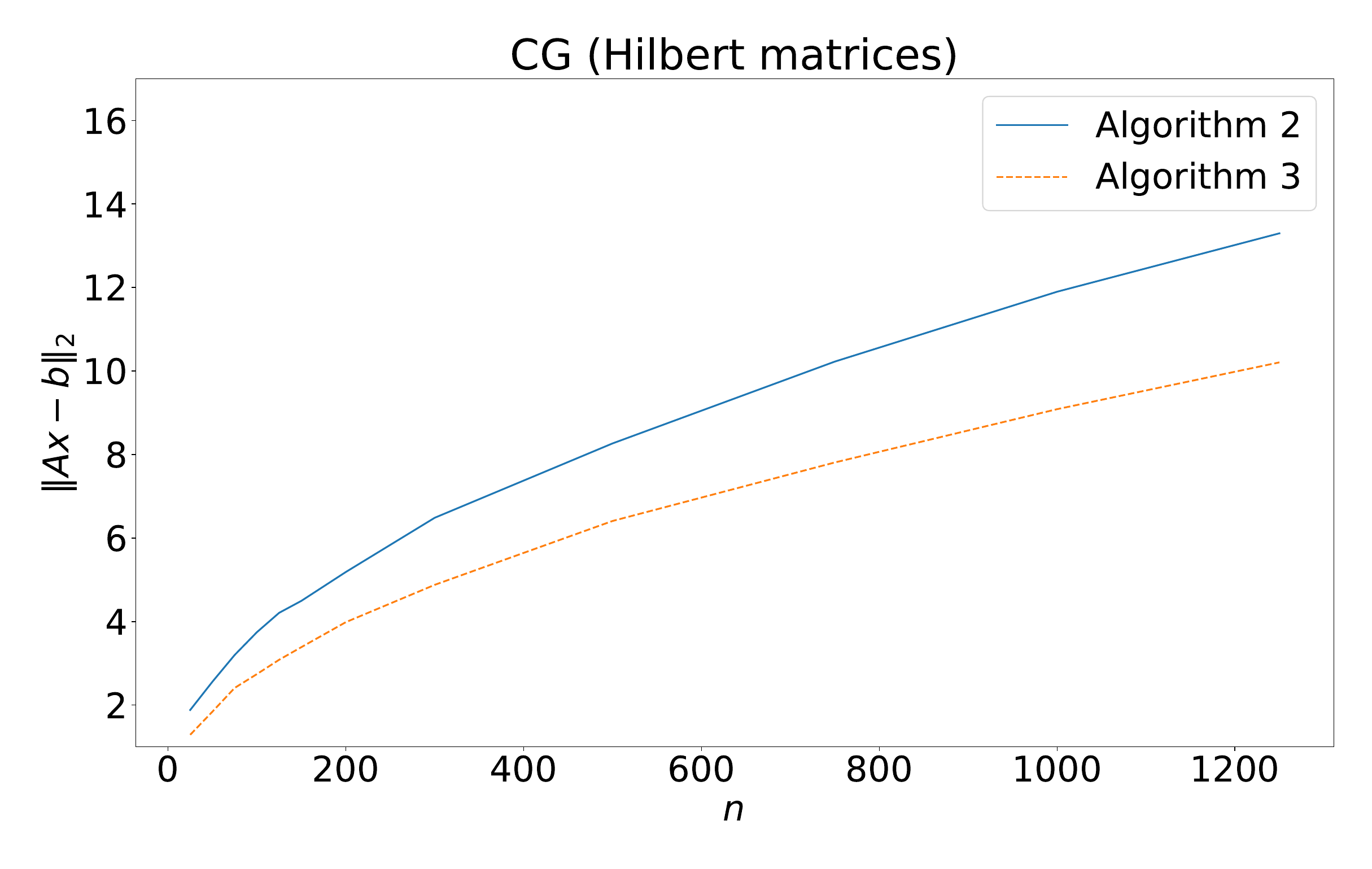}}
\caption{Residual norm of Algorithm~\ref{alg:two} vs.\ Algorithm~\ref{alg:three} with CG as the iterative solver for Hilbert matrices.}
\label{fig:CG2v3}
\end{figure}

\begin{figure}[htbp]
\centering
{\includegraphics[width=0.5\textwidth]{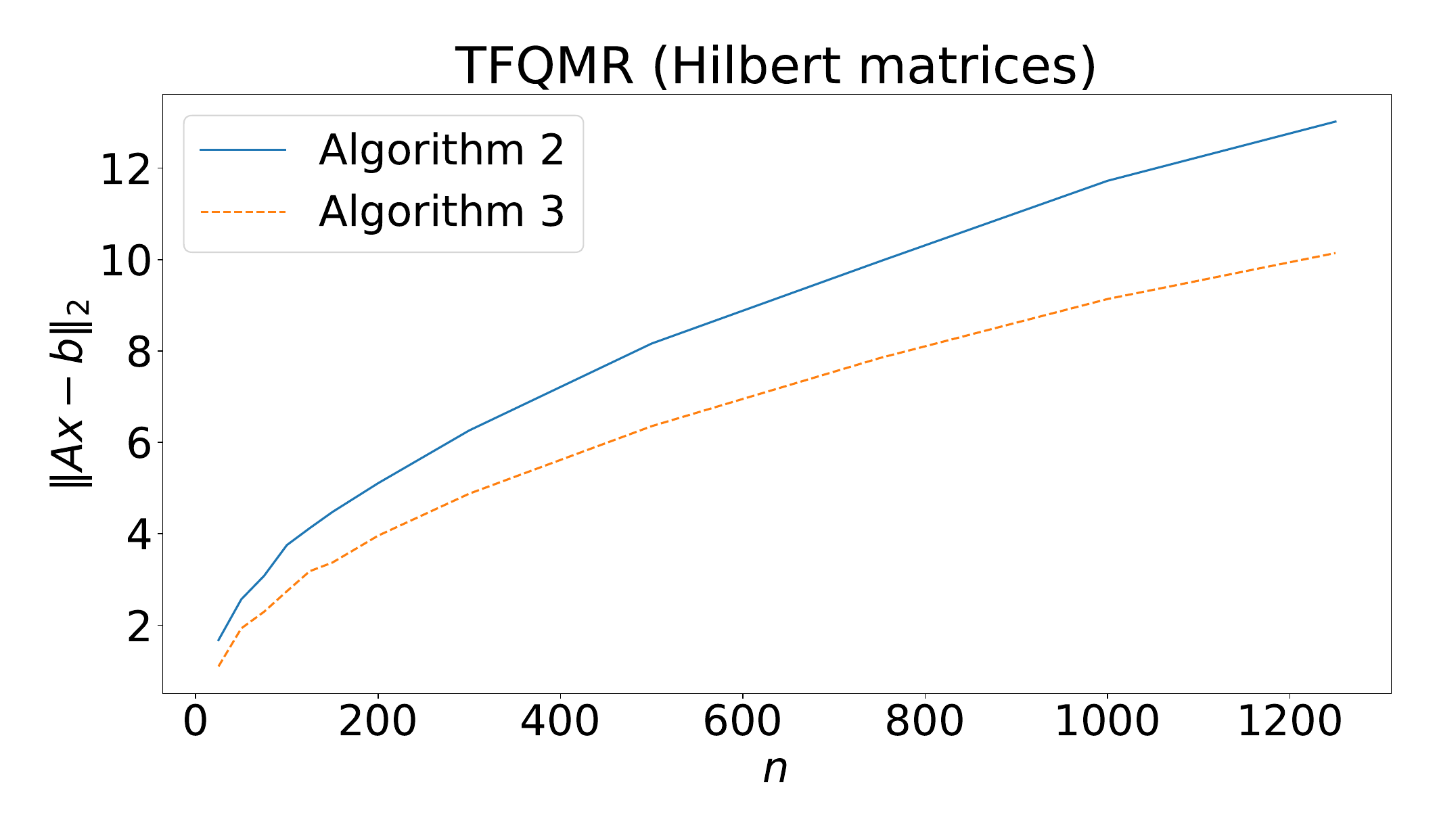}}
\caption{Residual norm of Algorithm~\ref{alg:two} vs.\ Algorithm~\ref{alg:three} with TFQMR as the iterative solver for Hilbert matrices.}
\label{fig:TFQMR2v3}
\end{figure}

Our second extended general stable iterative framework in Algorithm~\ref{alg:three} can provide further performance benefits over our first stable framework in Algorithm~\ref{alg:two} at the expense of a small additional computational cost, namely the solution of the $2\times 2$ system in Eq.~\eqref{eq:extended-framework}.
Representative samples of numerical experiments illustrating such additional performance benefits are presented in Figs.~\ref{fig:CG2v3}~--~\ref{fig:TFQMR2v3} for the CG and TFQMR algorithms respectively applied to Hilbert matrices, from which we observe that our extended stable Algorithm~\ref{alg:three} yields a smaller residual norm than our stable Algorithm~\ref{alg:two} where the gap in performance between the two algorithms grows as $n$ increases.

\subsection{Ill-conditioned random matrices}
We now turn to consider numerical experiments with respect to various ill-conditioned random matrices which allow us to explore certain behaviors of interest in a controlled manner.
Specifically, we construct such ill-conditioned matrices $\bA$ with a prescribed condition number $c$ by generating a random nonsingular symmetric matrix and rescaling\footnote{We consider the traditional condition number for the Euclidean norm, which is equal to the ratio of the extremal singular values.} the singular values to fall within the range $[1, c]$ \cite{kalantzis2018scalable}. 
In general,
our numerical results demonstrate that,
when the condition number is large, the classical iterative algorithms can diverge whereas the stable iterative algorithms based on our general frameworks do not diverge. 
Moreover, when the condition number is sufficiently small but not too small,
our numerical results show that the empirical convergence of the classical algorithms can be quite slow whereas the empirical convergence of our stable algorithms can be faster. 
These behaviors are readily apparent from across the large collection of numerical experiments performed within the context of random symmetric ill-conditioned matrices.
In particular, for many problem instances, 
the default GMRES algorithm terminates after exceeding the maximum number of iterations (typically $10n$ in \verb+SciPy+) without satisfying the standard stopping criteria, whereas our stable GMRES algorithm terminates after satisfying the stopping criteria well before reaching the the maximum number of iterations.

\begin{figure}[htbp]
\centering
{\includegraphics[width=0.5\textwidth]{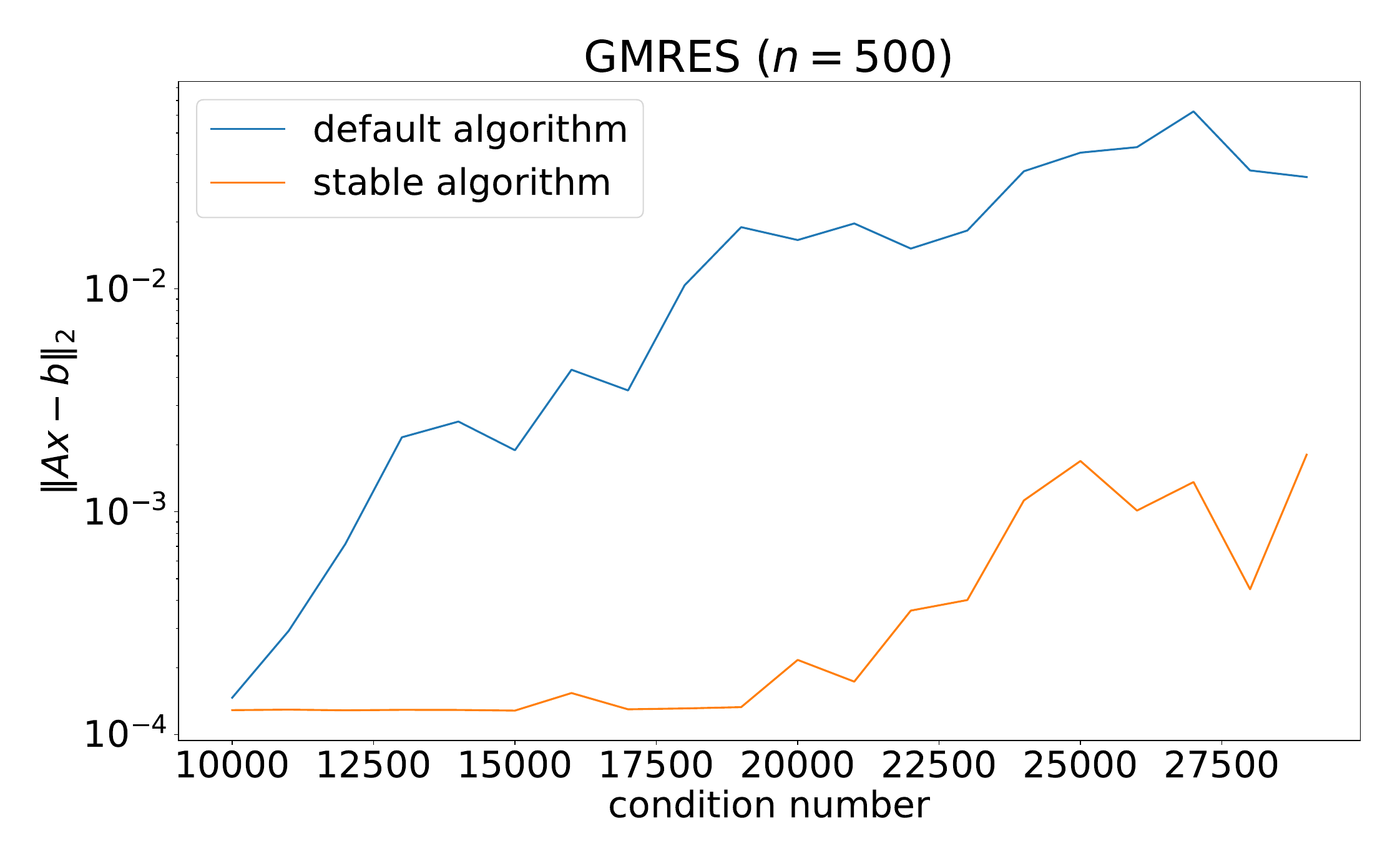}}
\caption{Solving random linear systems using GMRES where the stable algorithm is Algorithm \ref{alg:three}.}
\label{fig:gmres_cond}
\end{figure}

Fig.~\ref{fig:gmres_cond} shows a representative sample of the corresponding numerical results for the classical GMRES algorithm and the stable GMRES algorithm modified according to our second framework in Algorithm \ref{alg:three} with $n=500$ and varying condition number $c$.
We observe that, for such ill-conditioned systems, our stable GMRES algorithm can achieve a residual norm that is more than an order of magnitude better than the classical GMRES algorithm.

\begin{figure}[htbp]
\centering
{\includegraphics[width=0.5\textwidth]{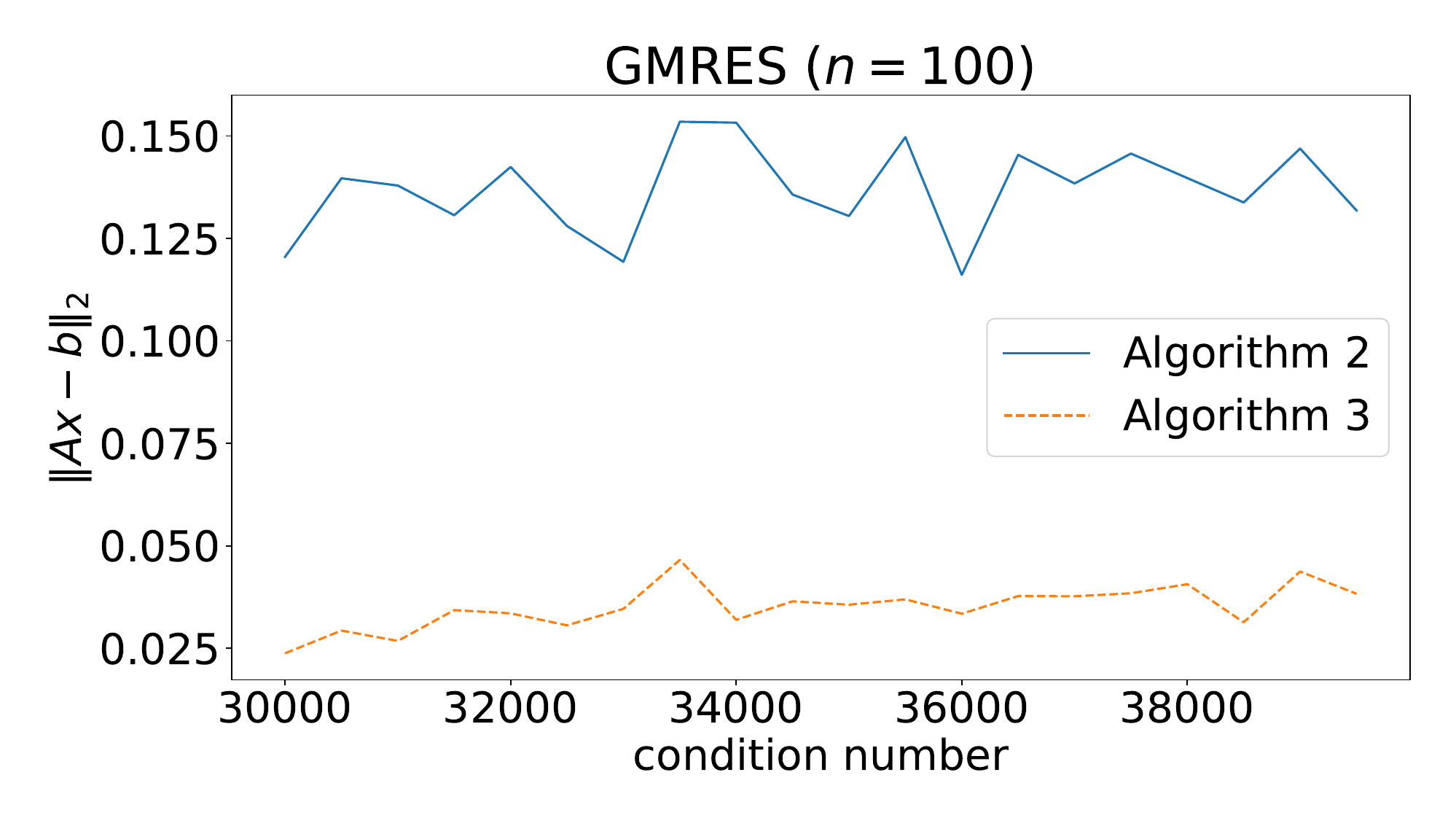}}
\caption{Residual norm of Algorithm~\ref{alg:two} vs.\ Algorithm~\ref{alg:three} with GMRES as the iterative solver for random matrices over $100$ trials.}
\label{fig:2v3}
\end{figure}

We again note that our second extended general stable iterative framework in Algorithm~\ref{alg:three} can provide further performance benefits over our first stable framework in Algorithm~\ref{alg:two} at the small expense of solving the $2\times 2$ system in Eq.~\eqref{eq:extended-framework}.
As another representative sample of numerical experiments illustrating such additional performance benefits, Fig.~\ref{fig:2v3} shows for the GMRES iterative solver applied to random matrices that our extended stable Algorithm~\ref{alg:three} yields a residual norm which is several times smaller than the residual norm under Algorithm~\ref{alg:two}.

\subsection{Real-world ill-conditioned matrices}
\begin{figure}[htbp]
\centering
{\includegraphics[width=0.5\textwidth]{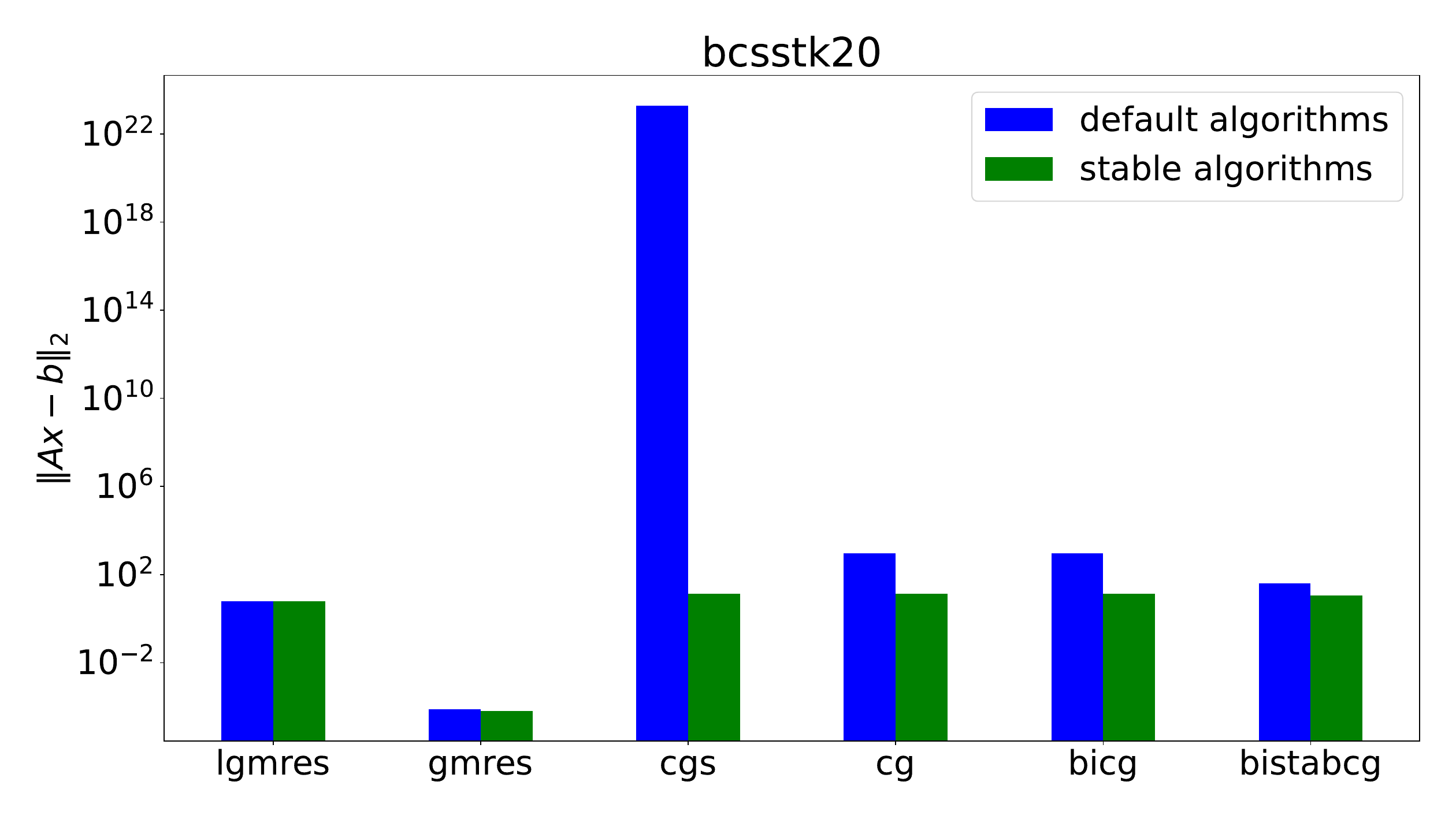}}
\caption{Solving matrix \texttt{bcsstk20}.}
\label{fig:bcsstk20}
\end{figure}

Lastly, we turn to study numerical experiments with respect to a set of real-world matrices taken from a standard library of matrices for real-world systems~\cite{sparse-ufl2011}
in which ill-conditioned linear systems are easily found.
To start, we first consider a particular real-world case where the ill-conditioned matrix $\bA$ is the symmetric positive definite sparse matrix \verb+bcsstk20+ of size $485\times 485$ obtained from a structural problem in the modeling of a suspension bridge frame. The condition number of this matrix $\bA$ is approximately $3.892662\times 10^{12}$.
A representative sample of our corresponding numerical experiments is presented in Fig.~\ref{fig:bcsstk20} which compares the residual norm of various default algorithms against the stable versions of these algorithms based on our first general framework.
We observe that our stable iterative algorithms perform equally well as the default algorithms \verb+lgmres+ and \verb+gmres+, and perform much better than the default algorithms \verb+cg+, \verb+bicg+, \verb+bicgstab+, and especially so for \verb+cgs+. 

\begin{figure}[htbp]
\centering
{\includegraphics[width=0.5\textwidth]{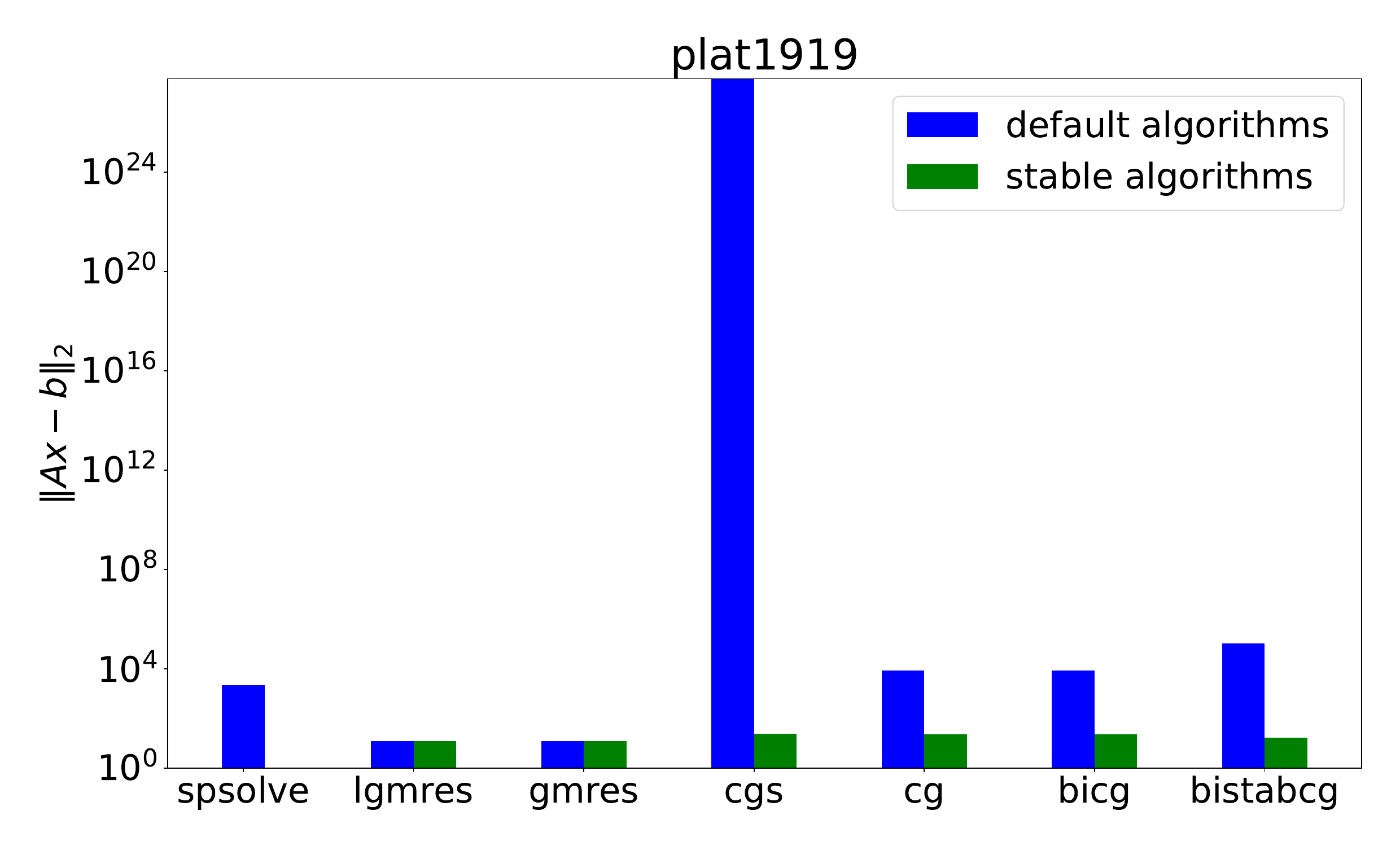}}
\caption{Solving matrix \texttt{plat1919}.}
\label{fig:plat1919}
\end{figure}

Next, we consider another particular real-world case where the ill-conditioned matrix $\bA$ is the singular symmetric sparse matrix \verb+plat1919+ of size $1919\times 1919$ obtained from a finite difference formulation of a three ocean problem. 
A representative sample of our corresponding numerical experiments is shown in Fig.~\ref{fig:plat1919} which compares the residual norm of various default algorithms against the stable versions of these algorithms based on our first general framework. 
Once again, we observe that our stable iterative algorithms perform equally well as the default algorithms \verb+lgmres+ and \verb+gmres+, and perform much better than the default algorithms \verb+spsolve+ (sparse direct solver), \verb+cg+, \verb+bicg+, \verb+bicgstab+, and especially so for \verb+cgs+. 

\section{Conclusions}
\label{sec:conclusions}
The serious issues of numerical overflow, inaccurate results and even divergence exhibited by classical Krylov subspace iterative solvers when applied to large-scale ill-conditioned linear systems represent a fundamental problem in high-performance computing that significantly reduces the applicability of these classical iterative methods in practice for obtaining solutions of such large-scale linear systems.
In this paper we mitigate this fundamental problem by introducing general algorithmic frameworks to augment the broad class of Krylov subspace iterative methods in a manner that solves large-scale ill-conditioned linear systems while ensuring stability and guaranteeing a lack of divergence. 
We apply our general algorithmic frameworks to current implementations of the iterative methods in \verb+SciPy+, and then conduct extensive numerical experiments that clearly demonstrate and quantify the significant extent to which classical Krylov subspace and related iterative methods diverge when the linear system is sufficiently ill-conditioned. 
Our numerical experiments across a wide range of synthetic and real-world ill-conditioned linear systems further show that the corresponding stable versions of these iterative methods augmented with our general algorithmic frameworks are stable, do not diverge, and provide significantly more accurate results. 
For less ill-conditioned linear system such that the classical Krylov iterative method does not diverge, our numerical results show that the Krylov subspace iterative methods augmented with our general algorithmic frameworks often provide significantly better accuracy than the original Krylov iterative method and never perform worse than the original method.

\end{document}